\documentstyle[11pt]{article}
\newtheorem{tm}{Theorem}[subsection]
\newtheorem{lm}[tm]{Lemma}
\newtheorem{pr}[tm]{Proposition}
\newtheorem{rmk}[tm]{Remark}
\newtheorem{cor}[tm]{Corollary}
\newtheorem{ex}[tm]{Example}

\newtheorem{??}[tm]{Question}
\newtheorem{defi}[tm]{Definition}

\font\tenmsb=msbm10
\font\sevenmsb=msbm7
\font\fivemsb=msbm5

\newfam\msbfam
\textfont\msbfam=\tenmsb
\scriptfont\msbfam=\sevenmsb
\scriptscriptfont\msbfam=\fivemsb
\def\Bbb#1{{\fam\msbfam #1}}

\font\teneufm=eufm10
\font\seveneufm=eufm7
\font\fiveeufm=eufm5
\newfam\eufmfam
\textfont\eufmfam=\teneufm
\scriptfont\eufmfam=\seveneufm
\scriptscriptfont\eufmfam=\fiveeufm

\def\lorw{\longrightarrow}

\newcommand\n{\noindent}

\newcommand{\nbigg}{\mathcal{G}}

\newcommand{\nbigk}{\mathcal{K}}

\newcommand{\nbigo}{\mathcal{O}}

\newcommand{\nbigs}{\mathcal{S}}

\newcommand{\seisuu}{\Bbb{Z}}

\newcommand{\lrarr}{\longrightarrow}

\def\ker{\mathop{\rm Ker}\nolimits}

\newcommand\ci{\cite}

\newcommand\rat{{\Bbb Q}}
\newcommand\comp{{\Bbb C}}

\newcommand\zed{{\Bbb Z}}

\newcommand\blacksquare{{\hspace*{\fill} $\Box$}}

\newcommand\Si{\Sigma}

\title{The Chow motive of semismall resolutions}
\author{
Mark Andrea A.  de Cataldo\thanks{
Partially supported by N.S.F. Grant DMS 9729992.}\, 
and Luca Migliorini\thanks{ member of GNSAGA, 
supported by MURST funds}
}

\date{March 15, 2002}

\begin{document}\maketitle
\abstract{
We consider proper, algebraic semismall maps
$f$
from a complex algebraic manifold $X.$ 
We show that the topological Decomposition Theorem
implies a ``motivic" decomposition theorem for the rational algebraic cycles
of $X$ and, in the case $X$ is compact, for  the Chow motive  
of $X.$     The result is a Chow-theoretic analogue
of Borho-MacPherson's observation \ci{b-m}
concerning the cohomology of the fibers
and their relation to the relevant strata for $f$.
Under suitable assumptions
on the stratification, we prove an explicit version of the
motivic decomposition theorem. The assumptions are fulfilled
in many cases of interest, e.g. in connection with resolutions 
of orbifolds and of some configuration spaces. We compute
the Chow motives and groups in some of these cases, e.g. 
the nested Hilbert schemes of points of a surface.
In an appendix with T. Mochizuki, we do the same for the 
parabolic Hilbert scheme of points on a surface.
The results above hold for mixed Hodge structures
and explain, in some cases, the equality
between orbifold Betti/Hodge numbers and  ordinary Betti/Hodge 
numbers for the crepant semismall resolutions in terms
of the existence of a natural map of mixed Hodge structures.
Most results hold over an algebraically closed field
and in the K\"ahler context.
}

\tableofcontents

\section{Introduction}\label{intro}
Borho and MacPherson \ci{b-m} first realized that the 
Decomposition Theorem
of Beilinson, Bernstein, Deligne and Gabber \ci{bbdg} 
has a simple and explicit statement
for complex proper algebraic {\em semismall } 
maps $f:X \to Y$ from a manifold. See 
$\S$\ref{semimap} for definitions and Theorem \ref{dtbbbd}.

Semismall maps appear naturally in a variety of contexts
arising in geometry and representation theory, e.g. the Springer 
resolutions of nilpotent cones, maps from  Gieseker to Uhlenbeck
moduli spaces, Hilbert-Chow maps for various Hilbert schemes
of points on surfaces, 
maps from Laumon to Drinfeld spaces of maps to flag
varieties, holomorphic symplectic resolutions, resolutions 
of theta divisors
for hyperelliptic curves, etc.

\n
Every stratum  for $f$ on $Y$  gives rise to
the local system given by the 
monodromy action on  the irreducible components of maximal dimension
of the fibers over the stratum.
Associated with this local system there is the corresponding
intersection cohomology complex
on the closure of the stratum.

\n
The Decomposition
Theorem for $f$ states that the direct image complex $Rf_* \rat_X[\dim{X}]$
is isomorphic to the direct sum 
over the set of
{\em relevant} strata
of these intersection cohomology complexes. 

\n 
A striking consequence
is that, for these maps, the cohomology of the fiber over
a point of a  stratum 
gets decomposed in, so to speak, the easiest possible way:
all but its  top cohomology is a direct sum
of the stalks of the cohomology sheaves
of the intersection cohomology complexes  of  the strata ``near" the point; 
the top cohomology comes, naturally,  from the stratum itself. 

\n
This statement is, in a sense, surprising. It is {\em as if}
the geometry of the fiber did not really count, but only 
the irreducible components of the fiber
and their monodromy along the stratum.

It is instructive to check how all the above fails
for the non-holomorphic map that contracts to a point the zero section
in the total space of  the trivial line bundle
on an algebraic curve.

\medskip
In this paper we show, using the Decomposition
Theorem, that given a map $f:X \to Y$ as above,
the picture for the Chow groups and the Chow motive of $X$
is  analogous. Here, we freely use  the language and properties
of correspondences and projectors even if
$X$ is not proper, in view of the refined Gysin formalism; see \ci{fulton}
and \ci{decatmig1}.

The theory of motives, of various kinds, has been conceived by Grothendieck as 
a candidate for a universal cohomology theory
for nonsingular projective varieties.
The so-called ``standard conjectures" 
on algebraic cycles are what is needed to realize
 the above vision.

\n
This paper does not touch upon this conjectural part of 
the theory of motives.
Rather, it proves decomposition results in the well-defined category of
Chow motives where
one considers nonsingular proper varieties $X, Y$ and the operations induced 
on their associated objects by algebraic cycles in $X \times Y.$
Natural maps such as blowing-ups, projections for projective bundles
and  finite maps
induce direct sum decompositions in this category. 

\medskip
\n
We show that new decompositions
arise for semismall maps. They imply analogous ones for
Chow groups, singular cohomology and mixed Hodge structures.

\medskip
Every relevant stratum with the associated  local system
gives rise to a non trivial 
correspondence in $X \times X$ which is a projector.
 This correspondence is constructed with the monodromy datum of
the irreducible components of the fibers over the stratum.

\n
These projectors are orthogonal to each other.
This is the Chow-theoretic and motivic analogue of the fact that
each relevant stratum contributes a direct summand in the Decomposition Theorem.

\n
The sum of these projectors is the diagonal of $X$. 
This corresponds
to the fact that the direct summands above add up to 
the direct image complex $Rf_* \rat_X
[\dim{X}].$

The precise formulation is given in Theorem \ref{motdec} which one may call
the {\em Motivic Decomposition Theorem} for semismall maps.

In general, even though the formulation of the Decomposition Theorem
may be quite explicit, as it is in the case of semismall maps, it may still
be difficult in practice to compute the intersection cohomology complexes
involved. It is important to know that the decomposition occurs, but it
 may be difficult to
extract information from that fact alone. The geometry of the particular situation
is instrumental in making this tool one of the most effective methods
of computation in the homology of algebraic maps.

This situation presents itself in the case of algebraic cycles and motives.

\n
On the one hand, Theorem \ref{motdec} allows to
talk about intersection Chow motives and intersection Chow groups
in the context of semismall resolutions; see
$\S$\ref{icmgg}. Let us remark that Chow motives
 whose realization is an intersection cohomology group
have been also defined in a special but highly significant case, 
that of parabolic 
cohomology for the $k$-fold product of the universal elliptic curve,
 in \cite{scholl} by  A. Scholl.

\n
On the other hand, it may be difficult
to compute with these objects.

\medskip
The main result of this paper is that,
under  suitable hypotheses
on the strata, which are verified in many cases of interest, one can
recover the algebraic cycles and the motive of $X$ in terms of the
corresponding objects for the strata.
The precise statement is given by Theorem 
\ref{maniii}, whose proof is independent
of Theorem \ref{motdec}.

A strengthening of this result
is given by Theorem \ref{mammasantissima}; though it applies
to more general situations, e.g. 
we use it to prove
Theorem \ref{supesymm},  its formulation involves
a projector that in practice may be difficult to compute. 

Let us state  the following easy consequence
of Theorem \ref{maniii}. It is sufficient for many applications
and  may give the reader a better idea of the results of this
paper
and of its  computational consequences.

Let $f:X \to Y$ be a semismall map of complex projective varieties with $X$
nonsingular.
For every connected relevant stratum 
$Y_{a},$ 
let $t_a: =(1/2)(\dim{X} - \dim{Y_a})$,  $y_{a} \in Y_{a}$ and
consider the $\pi_1(Y_{a},y_{a})-$set of irreducible maximal 
dimensional components of $f^{-1}(y_{a})$.
Let $\nu_{a}:Z_{a} \to Y_{a}$
the corresponding not necessarily connected  covering.
Assume that every $Z_{a}$ has a projective compactification
 $\overline{Z_{a}}$
which is the disjoint union of 
quotients of  nonsingular varieties such that the map $\nu_{a}$ extends
to a finite map $\overline{\nu_{a}}:\overline{Z_{a}} \to \overline{Y_{a}}$.
We denote the Chow motive of a nonsingular (or quotient) proper
variety $T$ by $[T]$, and the $n$-th twist of it by
the Lefschetz motive by $[T](n);$ e.g.
$[{{\Bbb P}^1}] \simeq [pt] \oplus [pt](1).$ 
We construct a correspondence $\overline{\Gamma} \subseteq
\coprod_a{ \overline{Z_a} \times X}$ and prove
\begin{tm}
\label{main}
There  are natural isomorphisms of Chow motives: 
$$
\overline{\Gamma}: 
\oplus_a \left[\, \overline{Z_{a}}\, \right] ( t_a) \simeq
[X] .
$$
\end{tm}
The analogous statement for Chow groups
follows easily.
The one for  Hodge structures, though  not new, as it  follows 
from Saito's \ci{sa},
is given here a direct proof
using  correspondences, i.e. without relying on the theory of
mixed Hodge modules.

\medskip
We provide the following applications of Theorem \ref{maniii} and Theorem
\ref{mammasantissima}.

\n
For every nonsingular complex surface,
we determine the Chow groups, the mixed Hodge structures and the Chow motive of
Hilbert schemes (Theorem \ref{cmhs2}), of nested Hilbert schemes
(Theorem \ref{cmnhs}) and, in an appendix
with T. Mochizuki ($\S$\ref{parhilb}), of parabolic Hilbert schemes of points.
All these results hold over any field and in the analytic
context. The former is not new; see \ci{decatmig1}
and also \ci{gogro}. The second refines G\"ottsche's result 
\ci{gogro}. The third is new.

\n
Our results allow to verify that for some semismall and crepant  resolutions of 
orbifolds, the mixed Hodge structure of the resolution is canonically isomorphic
to the corresponding orbifold structure. See $\S$\ref{wprdp}. This refines
the known fact that one had an equality
for the corresponding Hodge numbers. See \ci{donagi} and \ci{wangorb}.

\n
Another interesting
case, which we plan to consider in the near future,   in which the 
methods of this paper seem to be relevant,
is when a group acts by symplectic reflections on a 
nonsingular quasi projective variety
and the quotient is given
a symplectic resolution, cfr. Kaledin \ci{kal}, Verbitski \ci{verb}.

\bigskip
\n
{\bf Acknowledgments.}
We thank W. Wang for useful conversations.
 We thank A. Corti for stimulating conversations on the subject
 and for explaining to us the content of the paper \cite{ch}.
In a certain sense the main goal of this paper is to show that a small part of 
the mostly conjectural theory presented in \cite{ch}
works without relying on the standard conjectures. 
The first author thanks the Institute for Advanced Study, Princeton,
and the Mathematical Sciences Research Institute, Berkeley, 
for  their kind hospitality while this paper was being written.
T. Mochizuki, who has co-authored the appendix $\S$\ref{parhilb}
wishes to thank the Institute for Advanced Study, Princeton,
and to acknowledge 
N.S.F. Grant DMS 9729992.

\section{A motivic decomposition theorem for semismall maps}
\label{semimap}
We work with algebraic varieties over the 
field of complex numbers 
and consider Chow groups with 
rational coefficients.
The results of this paper hold, with minor changes, left to the reader,
over an algebraically closed field $K$ and in the complex analytic setting when
the map $f:X \to Y$ is K\"ahler, for the Decomposition Theorem holds
for $f$
in those contexts. See \ci{bbdg} and \ci{sa}; see  also \ci{decatmig2}
for a necessary and sufficient condition for
the  Decomposition Theorem to hold for $f$ in the analytic setting.
Over the field $K$ one must work with $\rat_l$-coefficients,
$l\neq char \, K,$ and with \'etale homology (cfr. \ci{Laumon}). However,
with  the exception of Theorem
\ref{motdec},  the relevant cycles have $\rat$-coefficients
and the results for cycles and motives remain valid as stated.
The results for (parabolic) Hilbert schemes 
hold over any field $F$ and in the complex analytic setting.

\subsection{Fibred products of semismall maps}
Throughout this section, $f: X \to Y$
is a proper surjective and {\em semismall} morphism
of complex algebraic varieties with $n:= \dim{X}.$ 

Recall that, in this context, the semismallness
of $f$
is equivalent to the following condition:
let $\delta \in {\Bbb N}$
and $Y^{\delta}_f:=
\{ y \in Y \, | \, \dim{f^{-1}(y)} = \delta \}$, then
$2\delta \leq \dim{Y} - \dim{Y^{\delta}_f},$ $ \forall \delta \geq 0.$

\n
On the other hand, the {\em smallness} of
$f$ is equivalent to requesting that
$2\delta < \dim{Y} - \dim{Y^{\delta}},$ $ \forall \delta > 0.$ 

By checking the condition on $\delta >0$, one sees immediately that
$f$ must be generically finite so that
$\dim{X}=\dim{Y}=n.$

The following elementary fact plays an important role
in what follows.
\begin{pr}
\label{prrr}
Let $f': X' \to Y$ be a proper surjective
and semismall map. Then
$$
\dim{X' \times_Y X }=n.
$$
\end{pr}
{\bf Proof.}
Let $\delta$, $\delta'$ $\in {\Bbb N}$
and define the locally closed subspaces $Y^{\delta \delta'}_{ff'}: =
Y^{\delta}_f \cap Y^{\delta'}_{f'}$.
We have that 
$Y= \coprod_{\delta\geq 0, \delta' \geq 0}{Y^{\delta \delta'}_{
ff'}
}
$.
Note that 
$$
X' \times_Y X= \coprod_{\delta \geq 0, \delta' \geq 0}
{
f'^{-1}(Y^{\delta \delta'}_{ff'}) \times_{Y^{\delta \delta'}_{ff'}}
f^{-1}(Y^{\delta \delta'}_{ff'})
},
$$
where all the subspaces on the right hand side are locally closed
in the left hand side.

\n
We have
$
\dim{ 
f'^{-1}(Y^{\delta \delta'}_{ff'}) \times_{ Y^{\delta \delta'}_{ff'} }
f^{-1}(Y^{\delta \delta'}_{ff'})  
}
= \delta + \delta' + \dim{ Y^{\delta \delta'}_{ff'} }
\leq n  - 
\frac{ \dim{ Y^{\delta_f} } + \dim{ Y^{\delta'_{f'}} } }{2}
+ \dim{ Y^{\delta \delta'}_{ff'} }
\leq n.
$
\blacksquare

\begin{rmk}
\label{components}
{\rm
A map $f:X \to Y$ is semismall 
if and ony if the irreducible components of $X \times_Y X$ 
have dimension at most equal to $\dim{X}$.
Similarly, a map $f:X \to Y$ is small if and only if 
the $n$-dimensional irreducible components of $X \times_Y X$
dominate $Y.$
}
\end{rmk}

\begin{rmk}
\label{serve}
{\rm
An argument analogous to the proof
of Proposition \ref{prrr} 
 shows that if we merely assume that
$f$ and $f'$ are small over their image and
$\dim{X} = \dim{X'},$ then 
$\dim { X\times_Y X'} \leq \dim{X} = \dim{X'}$ and equality holds
if and only if $f(X)=f'(X').$ 
}
\end{rmk}

\bigskip
There exists a  {\em stratification for $f:X \to Y$}. By this we mean
that           
$Y= \coprod_{b\in B}{Y_b}$, where every space $Y_b$
is a locally closed and smooth subvariety of $Y$
and the induced maps $f_|: f^{-1}(Y_b) \to Y_b$ 
are locally topologically trivial over
$Y_b.$ See \ci{gomabook}.

\begin{defi}
\label{rel}
{\em 
Let $A:= \{ a \in B \, | \, 
2\dim{f^{-1}(y)}= \dim{Y} - \dim{Y_a}, \; 
\forall y \in Y_a \}.$
We call $A$ the set of {\em relevant} strata for $f$ .
}
\end{defi}
\begin{rmk}
\label{bla}
{\rm 
While there may be many different stratifications for $f,$
the set of subvarieties 
$\{ \overline{Y_a} \}_{a \in A}$ is uniquely determined
by $f$.
}
\end{rmk}

\begin{rmk}
\label{smalluniquerel}
{\em
A small map has only one relevant stratum, i.e. the dense one.
}
\end{rmk}

\subsection{The intersection form
on the local system associated with
a
relevant stratum}
\label{intform}
In this section we assume that $X$ is a quotient variety,
i.e. that $X\simeq X'/G$ is the quotient of a nonsingular variety
$X'$
by the action of a finite group $G$  of automorphisms of $X'.$

\n
The goal is to define certain intersection numbers which can be defined
when $X$ is nonsingular and therefore also when $X$ is a quotient variety.
For more details on what follows see
\ci{decatmig2} and \ci{gomabook}.

Let $S$ be a relevant stratum for $f$ and $s\in S$.
We have $2\dim{f^{-1}(s)}=n - \dim{S}$. Let $N$ be a contractible 
euclidean neighborhood
of $s$ and $N'= f^{-1}(N)$. 
Let $\{F_i\}$
be the set of irreducible components of maximal dimension
$(1/2)(n + \dim{S})$ of $f^{-1}(S \cap N)$, and 
let $\{f_i\}$
be the set of irreducible components of maximal dimension
$(1/2)(n - \dim{S})$ of $f^{-1}(s)$. We denote by $f_i \cdot F_j$
the rational valued refined 
intersection product defined
in $N';$ see \ci{fulton}. 
It is independent of the choice of $s \in S \cap N$
and is a monodromy invariant, i.e. it is $\pi(S, s)$-invariant.

\begin{rmk}
\label{poschar}
{\rm
The  integer-valued intersection form is 
defined similarly over any algebraically closed field,  
provided an \'etale neighborhood 
is used to 
trivialize the monodromy of the components
of the fibre. 
}
\end{rmk}
We denote Borel-Moore homology (cfr.  \ci{cg} for example) with rational coefficients
by $H^{BM}(-).$
As $s$ varies on $S$, the Borel-Moore classes $[f_i]  \in
H^{BM}_{n -\dim{S}}
(f^{-1}(s))$ form a local system $L_S^*$
 and the intersection pairing above forms
a pairing $L_S^* \otimes L_S^* \to \rat_S$ on this local system.
Since the monodromy is finite, $L_S^*$ is abstractly, 
i.e. not necessarily via the pairing above, 
isomorphic to its dual
which can be identified with $L_S:=(R^{n - \dim{S}}f_* \rat_X)_{|S}.$

\subsection{The Decomposition Theorem
and homological correspondences}
\label{dthc}
In the case of constant coefficients,
the Decomposition Theorem
of Beilinson, Bernstein, Deligne and Gabber has a  simple statement
for semismall maps from a quotient variety.
For a different  proof for projective 
semismall maps
from projective complex manifolds
see \ci{decatmig2}.

Given a variety $Z$ and a local system $L$
on a Zariski-dense open subset of its regular part $Z_{reg}$,
we denote by $IC_Z(L)$ the  Goresky, MacPherson, Deligne 
intersection cohomology
complex associated with $L$. We denote $IC_Z(\rat_{Z_{reg}})$ simply by
$IC_Z$.

\begin{rmk}
\label{ns}
{\rm
If $Z$ is a  quotient variety, then $IC_Z\simeq \rat_Z[\dim{Z}].$
}
\end{rmk}

The following is a special case of the Decomposition Theorem
of Beilinson, Bernstein, Deligne and Gabber. 
See \ci{bbdg}, \ci{sa}, \ci{b-m}, \ci{decatmig2}. 
\begin{tm}
\label{dtbbbd}
Let $f: X \to Y$ be a proper, semismall algebraic map
of algebraic varieties with $X$ a quotient variety.
There is a canonical isomorphism
$$
Rf_*\rat_X[n] \simeq_{D^b_{cc}(Y)} \bigoplus_{a \in A} IC_{\overline{Y_a}}
(L_{a})
$$ 
where $L_{a}$ are the semisimple local systems on $Y_a$
given by the monodromy action on the maximal dimensional
irreducible components of the fibers of $f$ over $Y_a,$ and $D^b_{cc}(Y)$
is the bounded derived category of constructible complexes of 
sheaves of rational vector spaces on $Y.$
\end{tm}


\begin{lm}
\label{keylemma}
Let $g': U'\to V,$ $g:U \to V,$ $g'':U'' \to V$ be three proper 
maps
of analytic varieties with $U',$ $U$ and $U''$ 
quotients
of smooth varieties by the actions of   finite groups, and $i,$ $j$ and $k$
be integers.
There are  canonical isomorphisms
$$
\varphi_{U'U}: Hom_{D^b_{cc}(V)} (Rg'_* \rat_{U'} [i], 
Rg_* \rat_U [j]) \simeq H^{BM}_{2\dim{V} + i - j}(U' \times_V U).
$$
Given  morphisms $u: Rg'_* \rat_{U'} [i] \to 
Rg_* \rat_U [j]$, $v: Rg_* \rat_U [j]\to
Rg''_* \rat_{U''} [k]$, one has $\varphi_{U'U''} (v\circ u) =
\varphi_{UU''} (v) \bullet \varphi_{U'U} (u)$, 
where ``$\bullet$" denotes the refined 
composition
of homological correspondences.
\end{lm}
{\em Proof.} \ci{ch}, Lemma 2.21, Lemma 2.23.
\blacksquare

\subsection{A  motivic decomposition theorem
for semismall maps}
\label{mdtsm}
Let $X$ be a nonsingular algebraic variety of dimension $n$.
Given an algebraic variety $X',$
 a cycle $\gamma \in Z_l(X\times X')$, whose support
is proper over $X',$ defines, via the refined
Gysin formalism, a correspondence and a graded map of
cycles $\gamma_*: A_* (X) \to A_{* + l -n}(X').$  
This calculus is compatible with
composition, provided the appropriate ``properness'' conditions are imposed.
In addition, one may merely assume that $X$ is a quotient variety
by a finite group. See \ci{fulton}, \ci {decatmig1}.

We employ the standard notation for Chow motives; see
 \ci{manin} and \ci{fulton}. Later we shall switch to a more compact one.

The following is a formal consequence of Theorem \ref{dtbbbd} and Lemma
\ref{keylemma}.
\begin{tm}
\label{motdec}
Let $f: X \to Y$ be a proper, semismall, algebraic map
of complex algebraic varieties with $X$ a quotient variety
and $A$ be the set of relevant strata for $f.$
There exist algebraic cycles with
$\rat$-coefficients $C_a \in Z_n(X\times_Y X)$, $a \in A$,
such that

\n
i) $C_a\neq 0,$ $\forall a \in A,$

\n
ii) $\sum_a C_a = \Delta_X,$ in $Z_n(X\times_Y X),$

\n
iii) $C_a \circ C_{a'}= \delta_{aa'}C_a,$
as refined correspondences.

\n
In particular, there are canonical isomorphisms
of Chow groups
$$
A_*(X)= \bigoplus_{a\in A}{ {C_a}_* ( A(X) ) }
$$
and, when $X$ is proper, of Chow motives:
$$
(X, \Delta_X) \simeq \bigoplus_{a \in A}{(X,C_a)}.
$$
The cycles $C_a$ are supported on
$f^{-1} (\overline{Y_a})\times_{\overline{Y_a}}
f^{-1} (\overline{Y_a}).$ 
\end{tm}
{\em Proof.}
By Proposition \ref{prrr}
the class
map $cl: Z_n(X \times_Y X) \to H^{BM}_{2n}(X\times_Y X)$ is an isomorphism.
By Lemma \ref{keylemma} we get isomorphisms of 
$\rat$-algebras:
$$
End_{D^b_{cc}(Y)}(Rf_* \rat_{X}[n])\stackrel{\varphi}\longrightarrow
H^{BM}_{2n}(X\times_Y X)\stackrel{cl^{-1}}\longrightarrow
Z_n(X \times_Y X),
$$
where
 the algebra structures are given by composition
in  $D^b_{cc}(Y)$ for the first term and by 
the topological and algebraic/analytic Gysin formalism of correspondences
(which are compatible via the class map) for the second and third.

\n
Let $\phi: = cl^{-1} \circ \varphi.$

\n
By Theorem \ref{dtbbbd}
we have, $\forall a \in A$, canonical maps
$p_a: Rf_*\rat_X [n] \to IC_{\overline{Y_a}}(L_{a})$
considered as  elements of 
$End_{D^b_{cc}(Y)}(Rf_* \rat_{X}[n]).$

\n
Define $C_a : = \phi (p_a)$, $\forall a \in A.$

\n
The  three properties i), ii) and iii) for the  $C_a$
 follow from the analogous 
obvious properties of the $p_a$'s and by the compatibility
of $\varphi$ with compositions stated in Lemma \ref{keylemma}.

\n
The isomorphisms of Chow groups and Chow motives follow formally.

\n
The statement on the supports follows from the fact that the isomorphism
of Lemma \ref{keylemma} is compatible with restricting to open subsets of
 $Y$. If we take as an open set the complement of $\overline{Y_a}$, then
$p_a$ restricts to the zero map so that $C_a$ restricts to the zero cycle.
\blacksquare


\subsection{(Symmetric) products of semismall maps}
\label{spsm}

Let $f: X \to Y$ be semismall
with $X$ a quotient variety. It is easy to check that 
the induced maps
$f^l: X^l \to Y^l,$ $f^{(l)}: X^{(l)} \to Y^{(l)}$ 
and $f^{(\nu)} : X^{(\nu ) } \to Y^{(\nu ) }$ 
are semismall.

Let $A$ be the set of relevant strata for $f$. Then the sets
$A^l$, $A^{(l)}$ and $A^{(\nu)}$ are naturally labeling the
relevant strata for $f^l$, $f^{(l)}$ and $f^{(\nu )}.$

One has $Rf_* \rat_X [n]\simeq \oplus_{a\in A}IC_{\overline{Y_a}}(L_a)$
and $Rf^l_* \rat_{X^l}[ln] \simeq \oplus_{ \underline{a} \in A^l }
IC_{ \overline{ Y_{\underline{a} } } } (L_{\underline{a}}),$ etc.

Let $P \in Z_n(X\times X)$ have support that maps properly onto
either factor. Note that we are not assuming that $X$ is proper,
thus the assumption on the support. 
Let $P^2=P$ as a refined correspondence, i.e.
$P$ is a projector on $X$.
It induces projectors $P^l$, $P^{(l)}$ and 
$P^{(\nu)}$ on $X^l$, $X^{(l)}$ and $X^{(\nu )}$. See \ci{manin}, $\S2$.

If $P_a$ denotes the projector in $Z_n ( X\times_Y X )$ corresponding
to $IC_{\overline{Y_a}}(L_a)$, then
$P_{\underline{a}}$ corresponds to 
$IC_{\overline{Y_{\underline{a}}}} (L_{\underline{a}}).$
Similarly, for $f^{(l)}$ and for $f^{(\nu)}.$

\bigskip
We leave to the reader the task of making explicit,
using the set of indices $A^l$, $A^{(l)}$ 
and $A^{(\nu)}$ the analogues of
Theorem \ref{motdec}, \ref{maniii} and \ref{mammasantissima}
for the maps $f^l$, $f^{(l)}$ and $f^{(\nu)}$.

\section{Intersection Chow motives and groups in the presence of semismall resolutions}
\label{icmgg}
What follows is a generalization 
of the well known construction of Mumford \ci{mumford}, which we now recall:
Let $X$ be a normal surface with a singular point $x$ 
and let $r:\tilde{X} \to X$ be a resolution. 
Denote by $E_i$ the exceptional divisors in 
$r^{-1}(x)$ and by $\Lambda$ the inverse of the intersection matrix
$(E_i,E_j)$. 
Given a curve $C$ on $X$ one has a unique lifting to a curve $\tilde{C}$ 
on $\tilde{X}$ 
with the condition that $\tilde {C} \cdot E_i=0$ 
for all exceptional curves $E_i$. In other words,
 one is
considering the cycle $P=\Delta_{\tilde{X}}- \sum \Lambda_{i,j}
E_i \times E_j $ in $\tilde{X}\times_X\tilde{X}$,
which is easily seen to be a projector. A direct computation shows that the Chow motive
$(X,P)$ has the intersection cohomology of $X$ as its Betti realization.
In a completely analogous way, whenever a singular variety $Y$ 
is given  a semismall resolution
$f:X \to Y$, it is  
possible to construct a projector 
$P \in A(X \times_Y X)$ whose Betti realization is the intersection cohomology 
of $Y$. 
This construction is meaningful because the Chow motive turns out
 to be independent of the semismall resolution chosen, as we prove in 
Proposition \ref{indepofres}.

One can therefore talk about intersection Chow motives
and groups, in the presence of semismall resolutions.

Let $f: X \to Y \leftarrow X' : f'$ be two proper surjective
semismall maps, with $X$ and $X' $ $n$-dimensional
proper
nonsingular algebraic varieties. 
Assume that there exists a Zariski-dense open subset $U$
of $Y$ and an isomorphism $g: f^{-1}(U) \to f'^{-1}(U)$ such that
$f_|= f'_| \circ g$ and  that $f^{-1}(U) \to U \leftarrow f'^{-1}(U)$
are topological coverings with associated isomorphic local systems
$L_f$ and $L_{f'}$.

The projections $p_f: Rf_* \rat_X[n] \to IC_Y(L_f)=IC_Y(L_{f'})
\leftarrow Rf'_*\rat_{X'}[n] : p_{f'}$
induce, as in the proof of Theorem
\ref{motdec}, algebraic cycles $C_f \in Z_n(X\times_Y X)$ and 
$C_{f'} \in Z_n(X'\times_Y X')$ which give Chow motives
$(X,C_f)$ and $(X', C_{f'})$, direct summands of
$(X, \Delta_X)$ and $(X',\Delta_X')$, respectively.

\begin{pr}
\label{indepofres}
{\rm 
There is a canonical isomorphism of Chow motives
$$
C_{f'f}: (X', D_{f'}) \to (X, D_f)
$$
defined by an algebraic cycle $C_{f'f} \in Z_n (X' \times_Y X)$.
}
\end{pr}
{\em Proof.}
Let $p_{f'f}: Rf'_*\rat_{X'}[n] \to Rf_*\rat_{X}[n]$
be the map obtained as the composition of the natural maps:
$Rf'_*\rat_{X'}[n] \to IC_Y(L_{f'}) = IC_Y(L_f)\to Rf_*\rat_{X}[n].$
Consider the isomorphism of $\rat$-vector spaces $\psi_{X'X}: Hom_{D^b_{cc}(Y)}
(Rf'_*\rat_{X'}[n],Rf_*\rat_{X}[n])\simeq H^{BM}_{2n}(X'\times_Y X)
\simeq Z_n(X'\times_Y X)$ (see Lemma \ref{keylemma}
and
 Proposition \ref{prrr}) and define
$C_{f'f}:= \psi (p_{f'f})$.

\n
Analogously, define $p_{ff'}$ and $C_{ff'}$.

\n
Clearly, $p_{ff'}\circ p_{f'f}= p_{f'}$ and
$p_{f'f}\circ p_{ff'}= p_f$.
Since $\psi$ is compatible with compositions, see Lemma \ref{keylemma},
the same equations hold for the corresponding cycles, i.e.
$D_{ff'}\circ D_{f'f}= D_{f'}$ and
$D_{f'f}\circ D_{ff'}= D_f$.
The result follows by noting that
$D_{f'}$ and ${D_f}$ induce the identity
on the motives $(X, D_{f'})$ and $(X,D_f)$, respectively.
\blacksquare

\begin{rmk}
\label{inutile2}
{\rm
Any two small resolutions of an algebraic variety
have isomorphic motives, cycles and mixed Hodge structures.
This  can be shown directly using
the graph of the corresponding birational map.
}
\end{rmk}

\begin{rmk}
\label{icmot}
{\em
If $X\to Y \leftarrow X'$ are two semismall resolutions
(in general there is none), then $(X, D_f)$ and $(X',D_{f'})$
are both motivic counterparts of the topological 
$IC_Y$. Proposition \ref{indepofres} states that
they are canonically isomorphic and therefore independent
of the semismall resolution.
}
\end{rmk}

\section{An explicit motivic decomposition theorem}
\label{proof}
The purpose of this section is to prove Theorem
\ref{maniii}
which, though  less general than 
Theorem \ref{mammasantissima},  is quite explicit and affords
many applications.

\bigskip
From now on we denote the Chow motive of a nonsingular (or quotient) proper 
variety $T$ by $\left[ T \right]$ and  $\left[ T \right](n)$ will 
denote the
$n$-th twist by the Lefschetz motive.

\bigskip
The set-up is as follows.
Let $f: X \to Y$ be a proper, surjective and semismall
map of complex algebraic varieties with $X$ a quotient variety
and $A$ be the set of relevant strata for $f$.
Let $a \in A$, $Y_a$  be the corresponding relevant stratum,
$t_a= (1/2)( \dim{X}- \dim{Y_a} ),$
$y_a \in Y_a$, $E_a$ be the $\pi_1(Y_a, y_a)$-set
given by the irreducible components of maximal dimension
of $f^{-1} (y_a)$ 
and  $\nu_a :  Z_a= \coprod_{i}{Z_{a,i}} \to Y_a$ 
be the not necessarily connected \'etale covering
associated with $E_a.$
For every $a,$ $i,$ 
let $\overline{Z_{a,i}}$ be a complex algebraic variety
containing $Z_{a,i}$ as a Zariski-dense open subset and such that
there is a  {\em proper} surjective map $\overline{\nu_{a,i}} : \overline{ Z_{a,i} }
\to \overline{Y_a}$ extending $\nu_{a,i}:= {\nu_{a}}_{| Z_{a,i}}$, i.e. such that
$\overline{\nu_{a,i}}_{| Z_{a,i}} = \nu_{a,i}.$

\begin{tm}
\label{maniii}
Assume that, for every $a,$ $i$,  $\overline{ Z_{a,i} }$ is a quotient variety and 
$\overline{ \nu_{a,i} }$ is small. 
Then there exists a correspondence $\overline{\Gamma} \in 
Z_*( \coprod_a{ \overline{ Z_a } \times X } )$ inducing
isomorphisms of Chow groups 
$$
\overline{\Gamma}_* = \bigoplus_{a \in A}{ A_* ( \overline{Z_a} ) }  
\lorw A_*(X) \qquad (with \;A_*(\overline{Z_a} ) \lorw A_{* +t_a} (X),
$$
of mixed Hodge structures
$$
\overline{\Gamma}_*^H = \bigoplus_{a \in A}{ H^* ( \overline{Z_a} ) (t_a)}  
\lorw H^*(X)
$$
and, if $X$ is proper, of Chow motives.
$$
\overline{\Gamma}_* = \bigoplus_{a \in A}  \left[ \,\,{\overline{Z_a}}
\,\,\right] (t_a) 
\lorw \left[ X \right].
$$
\end{tm}

We shall prove the theorem
in this section by constructing $\overline{\Gamma}$ and its inverse explicitely.

\begin{rmk}
\label{small}
{\rm
Note that
if $g: U  \to V$ is a small proper map of analytic varieties, then  $Rf_*IC_U
\simeq IC_V(L)$, where 
$L$ is the local system given by the monodromy action
on the  points of a  general fiber of $f.$
}
\end{rmk}

\begin{rmk}
\label{dtconv}
{\rm
The assumptions on the varieties 
$\overline{Z_{a,i}}$ and the maps
$\overline{ \nu_{a,i} }$ imply, by 
 \ref{ns} and \ref{small}, that the 
Decomposition Theorem \ref{dtbbbd} reads as follows:
$$
Rf_*\rat_X[n]\simeq
\bigoplus_{a \in A}  R\overline{\nu_{a}}_* 
\rat_{\overline{Z}_{a} } [dim Y_{a}].
$$
}
\end{rmk}

\subsection{Construction of correspondences}
\label{corr}
We now introduce the correspondences
that define the isomorphism of motives.  
We shall  construct a set of correspondences
$\Gamma_{a,i,I}$ associated with a relevant stratum
$Y_{a}$. Their closures will define the wanted isomorphisms.

In order to simplify the notation, we temporarily work on one 
stratum at the time.
Let $S$ be a relevant stratum, $s \in S$,
$E$ be the right $G:=\pi_1(S,s)$-set 
of maximal dimensional irreducible components of $f^{-1}(s)$. 
Denote by $\nu_E:S_E \to S$ 
the not necessarily connected  \'etale cover of $S$ corresponding to $E$. 

\n
By  assumption, $S_E$ is
a  Zariski-dense open subset of 
a disjoint union of 
global quotient varieties $\overline{S_E}$, 
and $\nu_E: S_E  \to S$
 extends to a small map 
$\overline{\nu_E}: \overline{S_E} \to \overline{S}\subseteq Y.$
Note that  $\overline{S_E}$
is of pure dimension $\dim{S}.$ 
The correspondences we will use are supported on the closures in 
$\overline{S_E}\times_Y X$
of $S_E\times_Y X$. We introduce notation in order to deal with 
their top-dimensional irreducible components.

\n
From now on we use the term ``component" to  refer to a top-dimensional one. 
This should create no confusion.

\n
Fix representatives $o_1,\cdots,o_r$ of the $G$-orbits of $E$; 
this corresponds to fixing
 base points $s_1,\cdots,s_r$
 in the various connected
 components 
$S_1,\cdots,S_r$ of $S_E$. Let $G_i$ be the stabilizer of $o_i$.
The following is a well-known and elementary fact on $G$-sets:
\begin{lm}
\label{irr}
The irreducible components
of $S_i\times_Y X$ are in 1-1 correspondence with the $G_i$ orbits of $E.$
\end{lm}

It follows that in order 
to specify an irreducible component $\Gamma_{\,i,\,I}$ of $S_i\times_Y X$
it suffices to specify a $G_i$-orbit $I=\{f_{i_1},\cdots,f_{i_l}\}$
of components of $f^{-1}(s)$.

{\rm
Let $\Gamma_i$ denote   the correspondence 
$\Gamma_{i,o_i}$  associated with the $G_i$ orbit $\{o_i\},$ $i=1, \ldots r.$.
}

\n
Note that the fiber over $s_i$ of the natural map
$\Gamma_i \to S_i$ is naturally 
identified with the corresponding irreducible component
 of $f^{-1}(s)$ and that $\Gamma_i$ is an irreducible variety.

\n
Let $\Lambda=(\Lambda_{i,j})_{i,j \in E}$ be the inverse  of 
the intersection matrix associated 
with the pair  $(S,s)$. Clearly, $\Lambda_{ig,jg}=\Lambda_{i,j},$ 
$\forall \, g\in G,$
for the intersection numbers are monodromy invariant.
From this follows that,
for every $i \in E,$ 
$\sum_j \Lambda_{i,j}f_j$ is a  rational  linear combination of 
$G_i$-orbits in the set of maximal dimensional components
of $f^{-1}(s)$, and thus defines, by \ref{irr},
a correspondence $\Gamma'_i \in Z(S_i \times_Y X).$
Note that in general its support is not irreducible.

Let us summarize what we have done. In general, the pre-image of the stratum
$S$ is not irreducible. In addition, each irreducible
component may fail to be irreducible over small euclidean
 neighborhoods of $s \in S$. Each irreducible component
of the pre-image of a stratum corresponds to a connected \'etale cover
$S_i$ and the pull-back $\Gamma_i$ of this irreducible component
to $S_i$ is irreducible locally in the euclidean topology over $S_i$
and in fact has irreducible fibers over $S_i.$  
The intersection forms being non-degenerate allows to define
 $\Gamma'_i.$

\bigskip
We need the following simple
\begin{lm}
\label{ifzerothenzero}
Let $T \in Z_{\dim{S}}(S_i\times_S S_j)$ be a correspondence,
$p_i:S_i \to S \leftarrow S_j : p_j$ be the natural projections.
Assume that $T_*: A_0( \{s_i \} ) \lorw A_0 ( p_j(p_i^{-1}(s_i)) )$
is the zero map. Then $T=0$.
\end{lm}
{\em Proof.} 
Let $E_i$ and $E_j$ be the corresponding $\pi_1(S,s)$-sets.
The nonsingular space $S_i\times_S S_j$ is of pure dimension
$\dim{S}$ and has one connected component for every
$G_i$-orbit inside $E_j$. Note that $T$ is a rational linear combination
of these components, $T =\sum{q_t T_t}$, and that it defines a map
$T_*$ as in the statement purely by set theory, i.e. without recourse
to the refined Gysin formalism. We have that $T_*(s_i) = \sum q_t \{ T_t
\cap p_i^{-1}(s_i) \}.$ Since the $G_i$-orbits in $E_2$ are disjoint,
the vanishing of the left hand side implies that $r_t=0,$ $\forall t$.
The statement follows.
\blacksquare

\begin{lm}
\label{slm}
$^t\Gamma'_j\circ\Gamma_i=\delta_{ij}\Delta_{{S_i}}.$
\end{lm}
{\em Proof.} 
The composition of the two correspondences
is represented by a cycle of dimension equal to $\dim{S}$ supported on
${S_i }\times_Y {S_j}$.
By Lemma \ref{ifzerothenzero}, it is enough to prove that
$( ^t\Gamma'_j \circ \Gamma_i )_*(s_i)=\delta_{ij} (s_i).$ 
With the notation introduced in $\S$\ref{intform},
${\Gamma_i}_*(s_i)=f_i$ and 
${^t\Gamma_j'}_*(f_i)=(\sum\Lambda_{j,k}F_k \cdot f_i)(s_i) =\delta_{ij}s_i$.
\blacksquare

\bigskip
In what follows we need  to emphasize the stratum with which 
the correspondences are associated so we will revert to  the 
complete notation: we have $Z_a= \coprod_i{Z_{a,i}},$ 
$\Gamma_{a,i}$ and  $\Gamma'_{a,i}.$ 

\bigskip
Define $\overline{\Gamma_{a,i}}$ to be the closure of
$\Gamma_{a,i} \subseteq Z_{a,i}\times_Y X$ in $\overline{Z_{a,i}}\times_Y X.$
Similarely for $\overline{\Gamma_{a,i}'}.$
Let $\overline{\Gamma_a}:= \coprod_i{ \overline{\Gamma_{a,i}} },$ 
$\, \overline{\Gamma}: =  \coprod_a{ \overline{\Gamma_{a}} },$
$\, \overline{\Gamma'_a}:= \coprod_i{ \overline{\Gamma'_{a,i}} },$ 
$\, \overline{\Gamma'}: =  \coprod_a{ \overline{\Gamma'_{a}} }.$

\begin{lm}
\label{ddelta}
$^t\overline{ \Gamma'_{b,j} } \circ \overline{ \Gamma_{a,i} }=
\delta_{ab} \delta_{ij} \Delta_{
\overline{Z_{a,i}} \times \overline{ Z_{a,i} } }.$
\end{lm}
{\em Proof.} Let $a=b$. 
For reasons of dimension of supports, 
see a similar argument in \cite{decatmig1},
the composition of the two correspondences
is represented by a cycle of dimension equal to $\dim{S}$ supported on
$\overline{S_i }\times_Y \overline{S_j}$.
By the smallness assumption, the only such cycles 
are linear combinations of the maximal dimensional irreducible components 
of the fiber product, each of which  is  dominant over both factors; see
Remark \ref{components}. 
We conclude by an argument identical to the one in the proof
of Lemma \ref{slm}.

\n
Let $a\neq b.$ 
The composition is represented by a cycle of dimension  $(1/2) (\dim{Z_a}
+\dim{Z_b})$  supported 
on $\overline{Z_{a,i}}\times_Y \overline{Z_{b,j}}$ which, 
by an argument analogous
to the proof of Lemma \ref{prrr}, 
has dimension at most equal to 
$\min \{\dim{Z_{a,i}}, \dim{Z_{b,j}} \}.$
If $\dim{Z_a} \neq \dim{Z_b}$, then
the composition is zero because it is supported on too small a set.
If $\dim{Z_a} = \dim{Z_b}$, but $a\neq b$, then 
Remark \ref{serve} shows again that the support is 
too small and therefore the composition is trivial.
\blacksquare

\subsection{ End of proof of Theorem \ref{maniii}.}
\label{eod}
By Lemma \ref{ddelta} we see that
$^t\overline{\Gamma'} \circ \overline{\Gamma}= \Delta_{\coprod_{a,i} \overline{Z_{a,i}}}.$
By 
Remark \ref{dtconv}, it follows formally that
the map induced on $Rf_* \rat_X [n]$ by the cycle
$\overline{\Gamma'} \circ ^t\!\overline{\Gamma}$, 
via Lemma \ref{keylemma}, is the identity.

\n
The cycle $\overline{\Gamma} \circ ^t\!\overline{\Gamma'}$ is supported on the 
$\dim{X}$-part of  $X\times_Y X$, so that, by Lemma
 \ref{keylemma}, this cycle is the diagonal $\Delta_X.$
\blacksquare


\begin{rmk}
\label{related}
{\rm 
The projectors appearing in  Theorem \ref{motdec})
satisfy
$C_a =  \overline{\Gamma_a} \circ ^t\!\overline{\Gamma'_a}.$
See also Remark \ref{cagen}. 
}
\end{rmk}

\section{A strengthening of Theorem \ref{maniii}}
\label{mainone}
We place ourselves in the same situation as in $\S$\ref{proof},
but we shall allow the maps $\overline{\nu_{a,i}}$ to be semismall.
We can prove a stronger result at the price
of not being able to write down the correspondence
in an explicit form due to the presence of the projectors
stemming from the semismall maps not being small.

\bigskip
Let $a \in A$, $L_a$
be the associated local system. Note that
if $L_{a,i}:={\nu_{a,i}}_* \rat_{Z_{a,i}},$
then $L_a= \oplus_i{ L_{a,i}}.$

If we assume that the maps
 $\overline{\nu_{a,i}}$ are {\em semismall}, 
then, corresponding to the canonical projection
$R \overline{ \nu_{a,i} }_* \rat_{ \overline{Z_{a,i}} }[\dim{Y_a}_a] \lorw
IC_{\overline{Y_a} } (L_{a,i} ),$ we get, via
Theorem \ref{dtbbbd} and Lemma \ref{keylemma}
projectors
$P_{a,i} \in Z_{\dim{Y_a}}
( \overline{Z_{a,i}} \times_Y \overline{Z_{a,i}} ).$
Let $P_a= \coprod_i P_{a,i}$ and $P=\coprod_a P_a$.

Note that, as a refined correspondence, 
$P$ is a projector in $Z_* ( \coprod_{a,i,b,j} 
\overline{Z_{a,i}} \times_Y \overline{Z_{b,j}} )$ with trivial components
in  $Z_* ( 
\overline{Z_{a,i}} \times_Y \overline{Z_{b,j}} )$ whenever 
$a\neq b$ and/or $i\neq j.$ Similarly for $P_a$.
\begin{tm}
\label{mammasantissima}
Assume that the maps $\overline{\nu_{a,i}}$ are {\em semismall}.
Then
$$
P \circ \, ^t\overline{\Gamma'} \circ \overline{\Gamma} \circ P = P
$$
and
$$
\overline\Gamma \circ P \circ \, ^t\overline{\Gamma'} =\Delta_X.
$$
\end{tm}
{\em Proof.}
Let $p, \, p_a, \, p_{a,i}, \,\gamma, \,\gamma_a,\, \gamma_{a,i},$
$^t\gamma', \,^t\gamma'_a, \,^t\gamma'_{a,i}$, 
$\overline{\gamma}, \,\overline{\gamma_a},\, \overline{\gamma_{a,i}},$
$\overline{^t\gamma'}, \,\overline{^t\gamma'_a}, \,\overline{^t\gamma'_{a,i}},$
be the maps in $D^b_{cc}(Y)$ corresponding, via
Lemma \ref{keylemma}, to  the correspondences
labeled by the corresponding capital letters.

\n
By Lemma \ref{keylemma} we need to show that $p \circ
\overline{^t\gamma'} \circ \overline{\gamma} \circ p =p$.

\n
It is enough to show that $p_{a,i} \circ
\overline{^t\gamma'_{a,i}} \circ 
\overline{\gamma_{a,i}} \circ p_{a,i} =p_{a,i}$ $\forall a,\,i$.

\n
Since the first and the last map are $p_{a,i}$, it is enough to check
that the restriction 
$\overline{  ^t\gamma'_{ai}   }\circ \overline{ \gamma_{ai} }$
to $IC_{\overline{Y_a}}(L_{ai}$ is the identity.
This is true over 
$Y \setminus (\overline{Y_a} \setminus
Y_a)$ by the equation 
$
^t\Gamma'_{b,j}\circ \Gamma_{a,i} = \delta_{ab} \delta_{ij} \Delta_{Z_{a,i}}
$ 
and it is still true over $Y$ by the semisimplicity of the perverse sheaf
$IC_{\overline{Y_a}}(L_{a,i})$ in the abelian category of perverse sheaves 
on $Y.$ 

\n 
The second statement is proved similarly. It is enough to show
that $\overline{\gamma_a}\circ p_a \circ \overline{^t\gamma'}= \pi_a$,
where $\pi_a: Rf_* \rat_X [n] \simeq  \oplus_b IC_{\overline{Y_b}}(L_b) 
\lorw Rf_* \rat_X [n] \simeq  \oplus_b IC_{\overline{Y_b}}(L_b)$ 
is the natural 
projection  onto $IC_{\overline{Y_a}} (L_a)$.

\n
By reasons of semisimplicity $\overline{^t\gamma'_a}$
annihilates $ \oplus_{b\neq a} IC_{\overline{Y_b}}(L_b)$.
It follows that $p\circ \overline{^t\gamma'_a}=
\overline{^t\gamma'_a}$ so that
$\overline{\gamma_a}\circ p_a \circ \overline{^t\gamma'}=
\overline{\gamma_a} \circ  \overline{^t\gamma'_a}$.
By reason of semisimplicty again,
$\overline{\gamma_a} \circ  \overline{^t\gamma'_a}$
maps $Rf_* \rat_X[n]$ into $IC_{\overline{Y_a}}(L_a)$
and in fact is a-priori non-trivial only on the direct
summand $IC_{\overline{Y_a}}(L_a)$. We have shown above that
$ \overline{^t\gamma'_a} \circ \overline{\gamma_a}$
is the identity, when viewed as a self map
on $IC_{\overline{Y_a}}$. The statement 
follows
\blacksquare

\bigskip
 
The following corollaries and remarks follow.

\begin{cor}
\label{isochgr}
$$
{\overline{\Gamma}}_* :  
\bigoplus_{a \in A}{ {P_a}_* ( A_* (\overline{Z_a} )     ) }
=\bigoplus_{a \in A, i}{ {P_{a,i}}_* ( A_* (\overline{Z_{a,i}} )     ) }
\lorw A_*(X)
$$
is an isomorphism with inverse $ (P \circ \, ^t\overline{\Gamma'})_*=
P_* \circ \, ^t\overline{\Gamma'}_*.$
$$
{\overline{\Gamma}}_*^H :  
\bigoplus_{a \in A}{ {P_a}_* ( H^* (\overline{Z_a} )     ) }
=\bigoplus_{a \in A, i}{ {P_{a,i}}_* ( H^* (\overline{Z_{a,i}} )     ) }
\lorw H^*(X)
$$
is an isomorphism of mixed Hodge structures
with inverse $ (P \circ \, ^t\overline{\Gamma'})_*^H=
P_*^H \circ \, {^t\overline{\Gamma'}}^H_*.$
Let $X$ be proper.
Then 
$$
\overline{\Gamma} \circ P : \bigoplus_{a \in A}{
(\overline{Z_a}, P_a) ( t_a )}
= 
\bigoplus_{a,i}{
(\overline{Z_{a,i}}, P_{a,i}) ( t_a )}
 \lorw \left[\, X \,\right] 
$$
is an isomorphism of Chow motives with inverse $P \circ \, 
^t\overline{\Gamma'}.$
\end{cor}

\begin{rmk}
\label{shifts}
{\rm 
While the maps ${P_a}_*$  preserve degrees,
${\Gamma_a}_*$ increases degrees by $t_a.$
Similarly for ${P_a^H}{\!}_*$ and for $\Gamma_*^H.$
}
\end{rmk}

\begin{rmk}
\label{cagen}
{\rm 
In the context of Theorem \ref{mammasantissima},
the projectors $C_a$ of Theorem \ref{motdec} satisfy
$C_a = \overline{\Gamma_a} \circ P_a \circ \, ^t\overline{\Gamma'_a}.$
}
\end{rmk}

\begin{rmk}
\label{recsm}
{\rm
If the map $\overline{\nu_{a,i}}$ is small, then $P_{a,i}$
is the diagonal and it can be omitted. One recovers Theorem
\ref{maniii} at once. 
}
\end{rmk}


\section{Maps induced by maps between surfaces}
\label{mbs}
Let $X$ be a nonsingular, connected complex algebraic surface, 
$Y$ be an algebraic surface
and $f:X \to Y$ be a proper surjective holomorphic map. 
Note that $f$ is automatically semismall.

\n
The map $f$ induces semismall
maps $f_n: X^{[n]} \to Y^{(n)}.$ We write down the the relevant strata for 
$f_n$ and prove Theorem \ref{supesymm}, an explicit
version of  Theorem \ref{mammasantissima} . An application is given
in Theorem \ref{ale}.

\bigskip
Let $n$ be a positive integer,
$X^{[n]}$ the Hilbert scheme of $n$-points of $X$ (cfr. \ci{nakbook},
\ci{decatmig0}),
$S_n$ be the $n$-th symmetric group and $X^{(n)}:= X^n/S_n$ 
be the n-th symmetric product of $X.$ There is a natural crepant
and semismall map
$\pi_n: X^{[n]} \to X^{(n)}$. Let $P(n)$ be the set of partitions
$\nu$ of $n,$ i.e. 
$\nu=\nu_1 \geq \ldots \geq \nu_l$,
where every $\nu_k$ is a positive integer
and $\sum{\nu_k}=n.$
The same partition can also be denoted
by $a= (a_1, \ldots , a_n)$, where $a_i:=$ {\em 
the number of times
$i$ appears in $\nu.$} Clearly, $\sum{ia_i}=n$ and the length
$l(a)=l(\nu):=l=\sum{a_i}.$

\n
The morphism $\pi_n$ admits a stratification
$X^{[n]}= \coprod_{\nu \in P(n)}{X^{[n]}_{\nu}} \to
\coprod_{\nu \in P(n)}{X^{(n)}_{\nu}}= X^{(n)},$
 where every
$X^{(n)}_{\nu}:= \{ x \in X^{(n)} \, | \,
x= \sum{\nu_k x_k}, \, x_k \in X, \, x_k\neq x_{k'}
\; \hbox{ for }\; k \neq k' \},$
is a connected locally closed smooth subvariety
of $X^{(n)}$ of dimension $2l(\nu),$
$X^{[n]}_{\nu} := (\pi^{-1}( X^{(n)}_{\nu} ))_{red}$
and the fibers of $\pi$ over $X^{(n)}_{\nu}$
are irreducible of dimension
$n-l(\nu).$

\n
Let $X^{(\nu)}:=\prod_{i=1}^n X^{(a_i)}.$
There is a natural map $\pi_{\nu}: X^{(\nu)}
\to X^{(n)}$ which factors through the normalization of
$\overline{X^{(n)}_{\nu}} \subseteq X^{(n)}.$

\bigskip

The composition
$
f_n : X^{[n]} \stackrel{\pi_n}\lorw  X^{(n)} 
\stackrel{f^{(n)}}\lorw Y^{(n)}.
$
is semismall.

A stratification for
$f$ induces one for $f^{(n)}$. Using the known one for $\pi_n$
in conjunction with the one for $f^{(n)}$ one gets one for $f_n$.
For our purposes, we need only the stratification by the dimension of 
the fibers which requires less indices. 

Let $V:= Y^0_f \subseteq Y,$ $U:= f^{-1}V \subseteq X$
and $Y^1_f := \{y_1, \ldots , y_N \}$. Note that
$V$ contains the unique dense stratum for $f$, which is relevant
for $f$  and that
all the $y_k$ are relevant strata, if $Y^1_f \neq \emptyset.$
There are no more relevant strata.

Let $D_k := (f^{-1}(y_k))_{red}$. $D_k$ is one dimensional, but not necessarily
 pure-dimensional. The one-dimensional part of $D_k$
is a not necessarily connected configuration of curves
$\{C_{kh} \}_{h\in H_k}$ on $X.$

Let $N \neq 0.$
Define, for every $h \geq 0,$ 
$Q_N(h)= \{ \underline{m}= ( m_1, ..., m_N ) \, | \;
m_k \in \zed_{\geq 0}, \; \sum{m_k}=h \}.$

\n
If $N=0$, then define $Q_0(0)=\{\star \}$ to be a set with one element
and $Q_0 (h) = \emptyset, \; \forall h >0.$
 
We have a decomposition
$$
Y^{(n)} = \coprod_{i=0}^n \coprod_{\nu \in P(i)} 
\coprod_{\underline{m} \in Q(n-i)}
V^{(i)}_{\nu} \times \star_{ \underline{m} }
$$
where $\star$ denotes a point and we identify
$V^{(i)}_{\nu} \times \star_{ \underline{m} }$
with the locally closed subset of points of type
$\sum{\nu_j v_j} + \sum_{k=1}^N{ m_ky_k}.$

We have that $f_n^{-1} ( V^{(i)}_{\nu} \times \star_{ \underline{m} }  )
= U^{[i]}_{\nu} \times \prod_{k=1}^N{ \pi_{m_k}^{-1}( D_k^{(m_k)}   )},$
where both sides are given the reduced structure.

The reduced fibers over the corresponding stratum for $f_n$ are the product
of a  fiber of $\pi_i$ over $X^{(i)}_{\nu}$
with the reduced configuration  $\prod_{k=1}^N{ \pi_{m_k}^{-1}( D_k^{(m_k)}   )}$.

\n
Let $\Sigma$ be a finite set
and
$M_d(\Sigma)$ be the set of degee $d$ monomials on the  set
of variables $\Sigma.$

\n
The first factor above  is irreducible and the second has irreducible 
components of top dimension labeled by 
$M_{\underline{m}}:=\prod_{k=1}^N{ M_{m_k}(H_k) }.$
The monodromy on these components is trivial.

Consider the natural map $X^{(\nu)}= X^{(\nu)}_{\underline{m}} \to
Y^{(n)}$ obtained by composing
$f^{(\nu)}$ with the assignment
$(\,  z^1_1, \ldots  , z^1_{a1} \, ; \,  z^2_1 , \ldots , z^2_{a_2} \, ; \,
\ldots  \, ; \,   z^i_1, \ldots , z^i_{a_i} \, ) \lorw 
\sum_{j=1}^i \sum_{t=1}^{a_j}{ jz^j_t} + \sum_{k=1}^N{m_ky_k}.$  

Each one of this will count as a relevant stratum
as many times as the numbers of irreducible components over it.
At this point we can apply
Theorem \ref{mammasantissima} to $f_n$
and obtain precise analogues of 
the results of $\S$\ref{mainone}.
\begin{tm}
\label{supesymm}
There are natural isomorphisms of Chow groups
$$
 \bigoplus_{i=0}^n \bigoplus_{\nu \in P(i)} 
\bigoplus_{\underline{m} \in Q(n-i)}{
P^{(\nu)}_* ( A_*(X^{(\nu)}_{\underline{m}} ) )^{\oplus |M_{\underline{m}}|}
}
  \simeq
A_*( X^{[n]}),
$$
of mixed Hodge structures
$$
\bigoplus_{i=0}^n \bigoplus_{\nu \in P(i)} 
\bigoplus_{\underline{m} \in Q(n-i)}{
P^{(\nu)}_* ( H^{BM}_*(X^{(\nu)}_{\underline{m}} ) )^{\oplus |M_{\underline{m}}|}
(2n - i - l(\nu))} 
 \simeq
H^{BM}_*( X^{[n]}),
$$
and, when $X$ is assumed to be proper, of Chow motives

$$
 \bigoplus_{i=0}^n \bigoplus_{\nu \in P(i)} 
\bigoplus_{\underline{m} \in Q(n-i)}{
(X^{(\nu)}_{\underline{m}}, P^{(\nu)})^{\oplus |M_{\underline{m}}|}
( 2n-i -l(\nu) )
}
  \lorw      
[ \, X^{[n]}\,].
$$
\end{tm}

\section{Applications of Theorem \ref{maniii}}

We now give a series of applications of Theorem 
\ref{maniii} to the computation of Chow motives and groups;
see $\S$\ref{hsps}, $\S$\ref{nestmot}, $\S$\ref{wprdp} and
$\S$\ref{parhilb}. 
It will be sufficient to
determine the relevant strata for the maps in question, 
the fibers over them, and the monodromy
on the components of maximal dimension. Except for 
the situations in $\S$\ref{wprdp} and
$\S$\ref{parhilb}, the fibers are irreducible.
This fact alone implies trivial 
 monodromy.
Theorem \ref{maniii} applies to give the Chow motive and groups of
Hilbert schemes and of nested Hilbert schemes. 
In the former case, that determination had been 
done in \ci{decatmig1} using the same correspondences, but 
additional knowledge was required, such as
the affine cellular structure of the fibers, and the fact
that the pre-image of each stratum maps locally trivially over the stratum
in the \'etale topology. By way of contrast, here one only needs
to know that the fibers have only one component of maximal dimension.
For nested Hilbert schemes, the result is new. 
The fibers of the maps in  $\S$\ref{wprdp}
are not irreducible, but have trivial monodromy.
Finally, the fibers of the Hilbert-Chow map
for parabolic Hilbert schemes (cfr. $\S$\ref{parhilb}) are not irreducible,
but they are over the relevant strata and that is enough.

\subsection{The Hilbert scheme of points on a surface}
\label{hsps}
The following result follows from  either
Theorem \ref{mammasantissima}, Theorem \ref{maniii}
or even Theorem \ref{supesymm} (applied to the map $Id: X \to X$).

\begin{tm}
\label{cmhs2}
Let $X$ be a nonsingular algebraic surface. There are 
natural isomorphisms of Chow groups
$$
\overline{\Gamma}_* :  \bigoplus_{\nu \in P(n)}{
A_* (X^{(\nu)})   \simeq   A_*(X^{[n]}) ,
}
$$
of mixed Hodge structures
$$
{\overline{\Gamma}}^H_* :
 \bigoplus_{\nu \in P(n)}{
H^* (X^{(\nu)})( n-l(\nu) ) 
\simeq H^*( X^{[n]} )
}
$$
and, if $X$ is assumed to be proper, of Chow motives
$$
\overline{\Gamma} : 
 \bigoplus_{\nu \in P(n)}
[\, X^{(\nu)} \,](n-l(\nu)) 
\simeq
[\, X^{[n]} \,].
$$
\end{tm}  
\begin{rmk}
\label{kummer}
{\rm Let $A$ be an abelian surface and $K_n$ be its n-th generalized Kummer 
variety (\ci{beau}).
An application of Theorem \ref{maniii} gives a refinement at the level of 
Chow motives of the formulae
for the Betti and Hodge numbers of $A \times K_n$ given in \ci{go-so}, 
Theorems 7 and 8. We omit the details.}
\end{rmk}

\subsection{The motive of nested Hilbert schemes of surfaces
$X^{[n,n+1]}.$}
\label{nestmot}
We apply Theorem \ref{maniii} to nested Hilbert schemes.
The results we obtain refine G\"ottsche's \ci{gogro}.
Let $X$ be a nonsingular quasi-projective surface,
$n$ be a positive integer
and define
the nested Hilbert scheme of points on $X$
$X^{[n,n+1]}$ to be the closed subscheme
$\{ (\eta, \beta) \, | \, \eta \subseteq \beta \} 
\subseteq X^{[n]} \times X^{[n+1]}$ given its reduced structure.
As is well-known (cfr. \ci{cheah}), it is connected and smooth
of dimension $2n+2$ and admits a semismall map 
$\pi_{n,n+1}: X^{[n,n+1]} \to  
X^{(n)} \times X$ which is semismall.
It is not crepant, e.g. $X^{[1,2]} \to X \times X$ is the blowing 
up of the diagonal.
This map admits a natural stratification
which we now describe.

Fix $a \in P(n)$. Define $I_{a}= \, \{ \, 0  \, \} \, \coprod \,
 \{\,  j \; | \; 
a_{j} \neq 0  \, \}.$ 
Clearly, $I_{a} \setminus \{ 0 \} \neq \emptyset$.

Define subvarieties $X_{a,j} \subseteq X^{(n)} \times X$ as follows
$$
X_{a,j}= \left\{ \begin{array}{cl}
\{ (\zeta , x) \, | \; \zeta \in X^{(n)}_{a}, \; x \not\in |\zeta | \}  
&  \mbox{ if $j= 0,$   }    \\   
\{ (\zeta , x) \, | \; \zeta \in X^{(n)}_{a}, \; x \in |\zeta |,
\; length_{x} (\zeta) =j
\}   
& \mbox{otherwise.}
\end{array}
\right. 
$$
The subvarieties $X_{a,j} \subseteq 
X^{(n)} \times X$ are irreducible, smooth, locally closed of 
 dimension $2l(a)+2$, if $j=0$ and $2l(a)$ otherwise.

We have the following stratification
$$
X^{(n)} \times X = \coprod_{a \in P(n)} \coprod_{j \in I_{a}}{   
X_{a,j}
}
$$
inducing, as in $\S$\ref{mbs}, a stratification for
the map $\pi_{n,n+1}$ for which the fibers are irreducible.
The fibers over the stratum $X_{a,j}$ have dimension
$n-l(a),$ if $j=0$ and $n-l(a)+1,$ if $j \neq 0$.
All strata are relevant.

Define, for every $a \in P(n)$ and for every $j \in I_a$,
$X^{(a,j)}: =X^{a,j}/\Si_{a,j}$.
There is a natural map
$X^{(a,j)} \to X^{(n)}\times X$ which factors through the normalization
of $\overline{X_{a,j}} \subseteq X^{(n)} \times X.$

We are in the position to apply Theorem \ref{maniii}
and prove
\begin{tm}
\label{cmnhs}
Let $X$ be a nonsingular algebraic surface.
There are    isomorphisms of Chow groups
$$
\overline{\Gamma}_* : 
 \bigoplus_{a \in P(n), j \in I_a}{
A_*(X^{(a,j)})} \simeq  A_* (X^{[n,n+1]}),  
$$
of mixed Hodge structures 
$$
\overline{\Gamma}_*^H :
\bigoplus_{a \in P(n), j \in I_a}{
H^*(X^{(a,j)})}(m(a,j))
\simeq
H^* (X^{[n,n+1]})  
$$
and, if $X$ is proper, of Chow motives
$$
\overline{\Gamma} :
\bigoplus_{a \in P(n), j \in I_a}
[\, X^{(a,j)} \, ]  (m(a,j)) 
\simeq 
[ X^{[n,n+1]} ]  ,
$$
where
$m(a,j)= n - l(a)$ if $j=0$ and $m(a,j)=n - l(a) + 1$
if $j\neq 0$.
\end{tm}

\begin{rmk}
{\rm 
\ci{gogro} computes the class of $X^{[n,n+1]}$ in the 
Grothendieck ring of motives, when $X$ is defined over a field of
characteristic zero. 
The result affords the algebraic cycles modulo some standard conjectures
on algebraic cycles. 
\ci{cheah} computed the   hodge numbers
using virtual Hodge polynomials.
\ci{go-so} computed the mixed Hodge structure
using Saito's mixed Hodge modules.
}
\end{rmk}

\subsection{Wreath products, rational double points and orbifolds}
\label{wprdp}
Let $G$ be a finite group, $n$ be an integer and $G_n$ be
the associated Wreath product, i.e. the semidirect product
of $G^n$ and $S_n$. For background on what follows
see \ci{wangale} from which
all the constructions below are taken.
Let $G \subseteq SL_2(\comp )$ be a finite group. 
Their classification is known.
Let $ \tau: \widetilde{\comp^2/G} \to {\comp}^2/G$ be the
minimal and crepant resolution of the corresponding simple singularity.
Let $D: = \tau^{-1}(o)$ where $o \in Y$ is the singular point.
With the notation of this section, $N=1$, $y_1=0$
 and the divisor $D$ is supported on a tree
of nonsingular rational curves: $|D|= \cup_{k \in H_1}
C_k.$

There is a commutative diagram of semismall and crepant birational
morphisms
$$
\begin{array}{ccc}
\widetilde{\comp^2/G}^{[n]} & \stackrel{\pi_n}\lorw  & 
\widetilde{\comp^2/G}^{(n)} \\
\downarrow \tau_n & & \downarrow \tau^{(n)} \\
\comp^{2n}/G_n & \simeq & (\comp^2/G)^{(n)}.
\end{array}
$$
Let $\comp_0^2/G: = \comp^2/G \setminus \{ o \}.$
The following is a natural stratification for $\tau_n$ for which 
every stratum is relevant:
$$
\tau_n= \coprod{ {\tau}_{n,i,\nu} } : \coprod_{i=0}^n \coprod_{\nu \in P(i)}
{ (\comp_0^2/G)^{[i]}_{\nu} \times
\pi_{n-i}^{-1} ( D^{(n-i)}   )    }
\lorw
\coprod_{i=0}^n \coprod_{\nu \in P(i)}
{ {  (\comp_0^2/G)^{(i)}_{\nu} } } \times 
\{(n-i)o   \}.
$$
The fibers  over each stratum
$(\comp_0^2/G)^{(i)}_{\nu}  $ 
are isomorphic to the product of the corresponding 
irreducible fibers
of $\pi_{i, \nu}$ with  $\pi^{-1}_{n-i} ( D^{(n-i) } ).$
The latter factor has dimension $n-i$ and its irreducible 
components of dimension
$n-i$ are naturally labeled by
the set $M_{n-i}(H_1)$ of monomials of degree $n-i$ in the set of variables
$H_1$.

Let $P(i) \ni \nu = (a_1, \ldots, a_i )$. The stratum
$(\comp_0^2)^{(i)}_{\nu} \subseteq$ 
$ \prod_{j=1}^i{
(\comp^{2}/G)^{(a_j)}
} =
$ $\prod_{j=1}^i{\comp^{2a_j}/G_{a_j}}$ $
= \comp^{2l(\nu)}/\prod_{j=1}^i{G_{a_j}}.$

There is the  natural finite map
$\comp^{2l(\nu)}/\prod_{j=1}^i{G_{a_j}} \to
(\comp^2/G)^{(n)}$
which factors through the normalization of the closure
of the stratum
$  (\comp_0^2)^{(i)}_{\nu}$.

Denote by $P= \Delta_{   \widetilde{\comp^2} } -
\sum_{h,k}{\Lambda_{h,k} C_h \times C_k }$ be the projector corresponding
to the projection $R \tau_* \rat_{ \widetilde{\comp^2/G} } [2] \to
IC_{\comp^2/G}\simeq \rat_{\comp^2/G}[2].$
Note that one can compute explicitly the projectors
$P^{(\nu)}$, for every partition $\nu$ of every integer.
We can apply Theorem \ref{maniii} and Theorem \ref{supesymm}
and prove
\begin{tm}
\label{ale}
There are  canonical isomorphisms of Chow groups
$$
\bigoplus_{i=0}^n
\bigoplus_{\nu \in P(i)}
\bigoplus_{M_{n-i}( G_* \setminus \{ 1^- \} ) }
{ A_* ( \comp^{2l(\nu)}/\prod_{j=1}^i{G_{a_j}}   ) } 
\simeq 
$$
$$
\bigoplus_{i=0}^n
\bigoplus_{\nu \in P(i)}
\bigoplus_{M_{n-i}( G_* \setminus \{ 1^- \} ) }
P^{(\nu)}_* { A_* (\widetilde{\comp^2/G}^{(\nu) } )  } 
\simeq
$$
$$
 A(\widetilde{\comp^2/G}^{[n]}) \simeq
\bigoplus_{\mu \in P(n)}{ A(\widetilde{\comp^2/G}^{(\mu ) } )}.
$$
and of mixed Hodge structures
$$
\bigoplus_{i=0}^n
\bigoplus_{\nu \in P(i)}
\bigoplus_{M_{n-i}( G_* \setminus \{ 1^- \} ) }
{ H^* ( \comp^{2l(\nu)}/\prod_{j=1}^i{G_{a_j}}   ) } 
\lorw H^{*}(\widetilde{\comp^2/G}^{[n]}) \simeq
\bigoplus_{\mu \in P(n)}{ H^{*}(\widetilde{\comp^2/G}^{(\mu ) } )}.
$$
\end{tm}

 Let $G_*$ be the set of conjugacy classes of $G$.
There is a well-known bijection $b : G_* \setminus \{ 1^- \} 
\to T_G$ (cfr. \ci{itoreid}; see also \ci{McKay}).

\n
 Using the bijection $b: G_* \setminus \{ 1^- \} \to T_G$
we get a natural set of bijections
$M_{n-i}( G_* \setminus \{ 1^- \}   ) \to 
M_{n-i}(T_G)$.

\n
The set of triplets $(i, \nu , x)$ with $0 \leq i \leq n$, 
$\nu \in P(i)$, $x \in M_{n-i}( G_* \setminus \{ 1^- \}   )$
is in natural bijection with the set ${G_n}_*$
of conjugacy classes of $G_n$.

\n 
It follows easily that one can identify canonically
the orbifold mixed Hodge structure of the orbifold
$\comp^{2n}/G_n$ with the mixed Hodge structure
of $\widetilde{\comp^2/G}^{[n]}.$
The following result refines slightly \ci{wangorb}.
\begin{cor}
\label{orbiale}
There is a canonical isomorphism of mixed Hodge structures
$$
(H^*( \comp^{2n}/G_n ))_{orb} \simeq
H^* ( \widetilde{\comp^2/G}^{[n]} ). 
$$
\end{cor}

\begin{rmk}
\label{orbhnetc}
{\rm 
More generally,
let $Y'$ be a nonsingular algebraic surface on which a finite group
acts as a finite group of automorphism such
that the only fixed points are isolated points on $Y'.$
Let $Y=Y'/G$ and $f: X \to Y$ be its minimal resolution,
which is automatically crepant. One can prove, using the method above and
the analysis of twisted sectors in \ci{wangorb}, 
that the orbifold mixed Hodge structure
of $Y'^n/G_n$ is isomorphic to
the mixed Hodge structure of $X^{[n]}$. 
This refines slightly results in \ci{donagi}
and \ci{wangorb}.
}
\end{rmk}

\section{Appendix: the Chow groups and the Chow motive
of parabolic Hilbert schemes; with T. Mochizuki}
\label{parhilb}
\centerline{\bf M. de Cataldo, L. Migliorini, T. Mochizuki}

\subsection{The parabolic Hilbert scheme of points on a surface}
Let $X$ be an irreducible nonsingular surface
defined 
over an algebraically closed field
and $D$ be
a nonsingular curve on   $X$.
Let $I$ be an ideal sheaf on  $X$.
Recall that a {\em parabolic structure} on  $I$ of depth $h$ at  $D$
is a filtration of the ${\nbigo}_D$-modules, denoted here using the induced 
successive 
quotients
\[
 I\otimes {\cal O}_D    \lorw  {\cal G}_{1}   \lorw{\cal G}_{2}\lorw\cdots
 \lorw{\cal G}_{h}\lorw{\cal G}_{h+1}=0.
\]
An ideal sheaf  of points with a parabolic structure is
called a parabolic ideal of points.
The moduli scheme of the ideal sheaves of points with
parabolic structure is called the {\em parabolic Hilbert scheme
of points}.
For simplicity,
we only consider the case that  ${\cal G}_1$ is torsion. 
There is no loss of generality,
see \ci{mochi}.

We put the ${\cal K}_i=\ker({\cal G}_i\lorw {\cal G}_{i+1})$.
For a parabolic ideal sheaf,
we denote the length of  $\nbigk_{\alpha}$ 
by $l_{\alpha}$ for $\alpha=1,\ldots,h$.
We denote the $h$-tuple
$(l_1,\ldots,l_h)$
by $l_{\ast}$.
The data $(n,h,l_{\ast})$
is called the type of the parabolic ideal sheaf $(I,\nbigg_{\ast})$.

We denote the parabolic Hilbert scheme of points
of $X$ of type $(n,h,l_{\ast})$ by
$Hilb(X,D;n,h,l_{\ast})$.
It is irreducible, nonsingular of dimension
$2n + \sum_{\alpha =1}^h{l_{\alpha}}$ (cfr. \ci{mochi}).
Note that if $l_{\alpha}=0$ for every $\alpha$, then
$Hilb(X,D;n,h,l_{\ast})$
is isomorphic to $X^{[n]}.$


\subsection{The corresponding Hilbert-Chow morphism}
For any sheaf $\nbigg$
of finite length,
we denote its support with multiplicity
by $[\nbigg]$,
i.e.,
if the length of  $\nbigg$ at $x$ is $l(x)$,
then  $[\nbigg]$ defined to be $\sum l(x)\!\cdot\! x$.
We denote  $X^{(n)}\times \prod _i D^{(l_i)}$
by $X^{(n)}\times D^{(l_{\ast})}$.
We have the natural morphism $F_{n,l_{\ast}}:$
of $Hilb(X,D; n,h,l_{\ast}) \to
X^{(n)}\times D^{(l_{\ast})}$
defined by the assignment
$(I,\nbigg_{\ast})\longmapsto
([\nbigo_X/I],[\nbigk_{\alpha}])$.
We call it the {\em Hilbert-Chow} morphism.
\ci{mochi}, Corollary 3.1 imply the following
\begin{lm}
The morphism $F_{(n,l_{\ast})}$ is semismall.
\end{lm}

\begin{rmk}
\label{notirr}
{\rm
The fibers of $F_{n,l_*}$ are not irreducible, in general.
However, the ones over relevant strata are; see Lemma \ref{irrokg}.
}
\end{rmk}

\subsection{Stratification of the morphism $F_{(n,l_{\ast})}$}
We denote the set of integers larger than $i$ by
$\seisuu_{\geq\,i}$.
We put $A=\prod_{\alpha=0}^h \seisuu_{\geq\, 0}$.
An element $v$ of $A$ is described as $v=(v_{0},v_{\ast})$,
where $v_{0}$ is an element of $\seisuu_{\geq\,0}$
and $v_{\ast}=(v_1,\ldots,v_h)$
is an element of $\prod_{\alpha=1}^h\seisuu_{\geq\,0}$.
We denote the set $A-\{(0,\ldots,0)\}$ by $A'$.

We denote the set
$\{\chi:A'\lrarr \seisuu_{\geq \,0}\,|\,\sum_v\chi(v)<\infty\}$
by $\nbigs(A')$.
We have the natural morphism
$\Phi:\nbigs(A')\lrarr A$ defined by
$\chi\longmapsto \sum_v \chi(v)\!\cdot\! v$.
For any element $u$ of $A$,
we put ${\nbigs}(A',u):=\Phi^{-1}(u)$.

For any element $\chi\in \nbigs(A')$,
we set:
\[
 X_{\chi}:=
 \prod_{v\in A', \, v_{*}= 0}
 X^{(\chi(v))}
\times
 \prod_{v\in A', \, v_{*}\neq 0}
 D^{(\chi(v))}.
\]
A point of $X_{\chi}$
is denoted  $\bigl(\sum_{k=1}^{\chi(v)} x_{v\,k}\,|\,v\in A'\bigr)$,
where $x_{v\,k}$ is a point of $X$ (resp. $D$)
if $v_{\ast}=0$ (resp. $v_{\ast}\neq 0$).
The dimension of  $X_{\chi}$ is
$2\sum_{v_{\ast}=0}\chi(v)+\sum_{v_{\ast}\neq 0}\chi(v) $.
Let $X_{0\,\chi}$ denote
the open set of $X_{\chi}$
determined by the condition that $x_{v\,k}\neq x_{v'\,k'}$
for any $(v,k)\neq (v',k')$.

For any $n\in \seisuu_{\geq\,0}$ and
$l_{\ast}\in \prod_{\alpha=1}^h\seisuu_{\geq\,0}$,
we have the element $l_{\ast}(n)=(n,l_{\ast})\in A$.
If $\chi$ is contained in  ${\nbigs}(A',l_{\ast}(n))$,
we have the finite morphism $G_{\chi}$
of $X_{\chi}$ to $X^{(n)}\times D^{(l_{\ast})}$
defined as follows:
\[
 \Bigl(\,\sum_{k=1}^{\chi(v)} x_{v\,k}\,\Bigl|\Bigr.\,v\in A'\,\Bigr)
\longmapsto
 \Bigl(\,
 \sum_v \sum_{k=1}^{\chi(v)}v_{\alpha}x_{v\,k},\,
\Bigl|\Bigr.\,\alpha=0,1,\ldots, h
 \,\Bigr).
\]
The restriction of $G_{\chi}$
to  $X_{0\,\chi}$ is an embedding.
We denote the image $G_{\chi}(X_{0\,\chi})$
also by $X_{0\,\chi}$.
Then $\{X_{0\,\chi}\,|\,\chi\in {\nbigs}(A',l_{\ast}(n))\}$
gives a stratification  
for the morphism $F_{(n,l_{\ast})}$.

Let $e_{\alpha}$ denote the element of $A'$ whose $\beta$-th component
is defined to be $1\,\,(\alpha=\beta)$
or $0\,\,(\alpha\neq \beta)$.
We put
$C:=\{m\!\cdot\! e_{0}\,|\,m\in \seisuu_{\geq\,1}\}\cup
 \{m\!\cdot\! e_{0}+e_{\alpha}
  \,|\,m\in\seisuu_{\geq\,0},\,\alpha=1,\ldots,h\}$.
The following two lemmata follow
from \ci{mochi}, Corollary 3.1.
\begin{lm}
For any point $P\in X_{0\,\chi}$,
the dimension of the fiber $F_{(n,l_{\ast})}^{-1}(P)$
is $(n-\sum_{v_{\ast}=0}\chi(v))$.
The codimension of the stratum $X_{0\,\chi}$
is
$2(n-\sum_{v_{\ast}=0}\chi(v))+\sum_{\alpha}l_{\alpha}
   -\sum_{v_{\ast}\neq 0}\chi(v)$.
A stratum
 $X_{0\,\chi}$ is a relevant stratum
if and only if $\chi(v)=0$ for any 
$v$ which is not contained in $C$.
\end{lm}
\begin{lm}
\label{irrokg}
If a  stratum $X_{0\,\chi}$ is relevant,
then the  top dimensional part of the inverse image
$F_{(n,l_{\ast})}^{-1}(X_{0\,\chi})$ is irreducible.
\end{lm}
Let $\bar{\nbigs}(A',l_{\ast}(n))$ denote the subset
of $\chi\in {\nbigs}(A',l_{\ast}(n))$ given by the
relevant strata.
The set $\bar{\nbigs}(A',l_{\ast}(n))$
can be regarded as the set of $\seisuu_{\geq\,0}$-valued
functions $\chi$ on  $C$
satisfying $\sum_{v\in C}\chi(v)v=l_{\ast}(n)$.

\begin{ex}
\label{exampleone}
{\rm 
($h=1,n=1,l_1=1$)\\
In this case, the parabolic Hilbert scheme is
isomorphic to the blowing up of  $X\times D$
with center $D$,
where  $D$ is embedded by the composition 
of the diagonal embedding $D\subset D\times D$ with
the natural inclusion $D\times D\subset X\times D$.
The Hilbert-Chow morphism is the blowing up morphism.
The set ${\nbigs}(A',(1,1))$ has two elements:
$\delta_{(0,1)}+\delta_{(1,0)}$ and $\delta_{(1,1)}$.
Here the function $\delta_v:A'\lrarr\seisuu_{\geq\,0}$
is defined by $\delta_v(w)=0\,\,(w\neq v)$ and $\delta_v(w)=1\,\,(w=v)$.
In this case $\bar{\nbigs}(A',(1,1))={\nbigs}(A',(1,1)).$
The stratum corresponding to  $\delta_{(1,0)}+\delta_{(0,1)}$
(resp. the $\delta_{(1,1)}$)
is $X\times D-D$ (resp. the $D$).
The dimensions of the fibers corresponding to
$\delta_{(1,0)}+\delta_{(0,1)}$
and $\delta_{(1,1)}$
are $0$ and $1,$ respectively.
}
\end{ex}
\begin{ex}
\label{exampletwo}
{\rm ($h=1,n=1,l_1=2$)\\
In this case, ${\nbigs}(A',(1,2))$ 
has four elements: $\chi_1:=\delta_{(1,0)}+2\delta_{(0,1)}$,
$\chi_2:=\delta_{(1,0)}+\delta_{(0,2)}$,
$\chi_3:=\delta_{(1,1)}+\delta_{(0,1)}$
and $\chi_4:=\delta_{(1,2)}$.
The dimensions $d_i$ of the fibers 
and the codimensions  $c_i$ of the strata
corresponding to  $\chi_i$ are as follows:
$(d_1,c_1)=(0,0),\,
 (d_2,c_2)=(0,1),\,
 (d_3,c_3)=(1,2),\,
 (d_4,c_4)=(1,3)$.
Thus 
 $\bar{\nbigs}(A',(1,2))$ has two elements
$\chi_1$ and $\chi_3$.
}
\end{ex}

\subsection{The results}
Let $\overline{\Gamma}$ be the correspondence constructed 
using the relevant strata above and $\S$\ref{proof}.
Let $t_{\chi}= n- \sum_{v_*=0} \chi (v)$.

\begin{tm}
\label{hilbpartuch}
We have natural isomorphisms of Chow groups
$$
\overline{\Gamma}_*: \bigoplus_{\chi \in \overline{S}(A', l_*(n))}{
A_*(X_{\chi}) \lorw A_* ( Hilb(X, D ; n,h,l_*) ) },
$$
of mixed Hodge structures
$$
\overline{\Gamma}_*^H: \bigoplus_{\chi \in \overline{S}(A', l_*(n))}{
H^*(X_{\chi}) (t_{\chi}) 
\lorw H^* ( Hilb(X, D ; n,h,l_*) ) },
$$
and, if $X$ is assumed to be proper, of Chow motives
$$
\overline{\Gamma}: \bigoplus_{\chi \in \overline{S}(A', l_*(n)) }{
[\, X_{\chi} \,]   (t_{\chi}) }  \lorw
[\, Hilb(X, D ; n,h,l_*) \,] .
$$
\end{tm}

In what follows, for an algebraic variety
$Y$,   $b_i(Y)$ denotes the $i$-th Betti number of $Y$,
$P(Y)(z)=\sum_i b_i(Y)z^i$ denotes the Poincar\'e polynomial of $Y,$
$h^{p,q}(Y)$ denote the $(p,q)$-th Hodge number of $Y$,
$h(Y)(x,y)=\sum_{x,y}h^{p,q}(Y)x^py^q$ denotes the Hodge polynomial of
$Y$,
and we put $\epsilon(p,q)=(-1)^{p+q-1}$.

\begin{cor}
\label{coronform}
The generating functions
for Chow motives (where we denote a motive
$(T, \Delta_T)$ simply by $[T]$), Betti and Hodge numbers are
$$
\sum_{n,l_{\ast}}
 [\, Hilb(X,D; n,h,l_{\ast}) \,]\cdot
 t^n\prod_{\alpha=1}^h s_{\alpha}^{l_{\alpha}}
= 
$$
$$
\prod_{i=1}^{\infty}
 \Bigl(\sum_{m=0}^{\infty} [\, X^{(m)}\,] ((i-1)m) \cdot t^{im}\Bigr)
\times
 \prod_{\alpha=1}^h
 \prod_{i=0}^{\infty}
 \Bigl(\sum_{m=0}^{\infty} [\, D^{(m)} ]  (im) \cdot t^{im}s_{\alpha}^m
 \Bigr),
$$
$$
\sum_{n,l_{\ast}}
 P(Hilb(X,D;n,h,l_{\ast}))(z)\cdot t^n
 \prod_{\alpha=1}^h s_{\alpha}^{l_{\alpha}}
=
$$
$$
\prod_{m\geq 1}
 \frac{(1+z^{2m-1}t^m)^{b_1(X)}(1+z^{2m+1}t^m)^{b_3(X)}}
 {(1-z^{2m-2}t^m)^{b_0(X)}(1-z^{2m}t^m)^{b_2(X)}(1-z^{2m+2}t^m)^{b_4(X)}}\\
\times
$$
$$
\prod_{\alpha=1}^h
\prod_{m\geq 0}
\frac{(1+z^{2m+1}t^{m}s_{\alpha})^{b_1(D)}}
      {(1-z^{2m}t^{m}s_{\alpha})^{b_0(D)}
   (1-z^{2m+2}t^{m}s_{\alpha})^{b_2(D)}}.
$$
$$
\sum_{n,l_{\ast}}
 h(Hilb(X,D,n,h,l_{\ast}))(x,y)
  t^n\prod_{\alpha=1}^hs_{\alpha}^{l_{\alpha}}
=
$$
$$
 \prod_{m\geq 1}
 \prod_{0 \leq p,q \leq 2}
  \Bigl(1+\epsilon(p,q)x^{p+m-1}y^{q+m-1}t^m
  \Bigr)^{\epsilon(p,q)\cdot h^{p,q}(X)}\\
\times
$$
$$
\prod_{\alpha=1}^h
\prod_{m\geq 0}
 \prod_{0 \leq p,q \leq 1}
  \Bigl(1+\epsilon(p,q)x^{p+m}y^{q+m}t^m s_{\alpha}
  \Bigr)^{\epsilon(p,q)\cdot h^{p,q}(D)}
$$
\end{cor}

\begin{rmk}
\label{alrea}
{\rm 
The generating functions for the Betti and Hodge numbers can be found in
\ci{mochi}. The slightly weaker statement in the Grothendieck
ring of motives can also be deduced easily from
\ci{mochi}.
}
\end{rmk}

\end{document}